\def\vers{Jan. 12, 2006, v.3}
\magnification=1200
\hsize=6.5truein
\vsize=8.9truein
\font\bigfont=cmr10 at 14pt
\font\mfont=cmr9
\font\sfont=cmr8
\font\sbfont=cmbx8
\font\sifont=cmti8
\def\scirc{\raise.2ex\hbox{${\scriptstyle\circ}$}}
\def\ssbull{\raise.2ex\hbox{${\scriptscriptstyle\bullet}$}}
\def\msum{\hbox{$\sum$}}
\def\mopls{\hbox{$\bigoplus$}}
\def\wdg{\hbox{$\wedge$}}
\def\bwdg{\hbox{$\bigwedge$}}
\def\bA{{\bf A}}
\def\bC{{\bf C}}
\def\be{{\bf e}}
\def\bN{{\bf N}}
\def\bP{{\bf P}}
\def\bQ{{\bf Q}}
\def\bZ{{\bf Z}}
\def\cA{{\cal A}}
\def\cB{{\cal B}}
\def\cD{{\cal D}}
\def\cG{{\cal G}}
\def\cH{{\cal H}}
\def\cI{{\cal I}}
\def\cJ{{\cal J}}
\def\cO{{\cal O}}
\def\cS{{\cal S}}

\def\oF{\overline{F}}

\def\ocH{\overline{\cal H}}
\def\ocS{\overline{\cal S}}

\def\tV{\widetilde{V}}
\def\tX{\widetilde{X}}
\def\tY{\widetilde{Y}}

\def\tom{\widetilde{\omega}}
\def\cSpec{{{\cal S}pec}}
\def\Im{\hbox{\rm Im}}
\def\Ker{\hbox{\rm Ker}}

\def\DR{\hbox{\rm DR}}
\def\Gr{\hbox{\rm Gr}}
\def\Sing{\hbox{\rm Sing}}
\def\rd{\partial}
\def\IH{\hbox{\rm IH}}
\def\simto{\buildrel\sim\over\longrightarrow}
\def\dfoverf{\hbox{${df\over f}$}}

\hbox{}
\vskip 1cm

\centerline{\bigfont A Generalization of Griffiths' Theorem}

\smallskip
\centerline{\bigfont on Rational Integrals}

\bigskip
\centerline{Alexandru Dimca and Morihiko Saito}

\bigskip\bigskip
{\narrower\noindent
{\sbfont Abstract.} {\sfont
We generalize Griffiths' theorem on the Hodge filtration of
the primitive cohomology of a smooth projective hypersurface,
using the local Bernstein-Sato polynomials, the
$ V $-filtration of Kashiwara and Malgrange along the hypersurface
and the Brieskorn module of the global defining equation of the
hypersurface.}
\par}

\bigskip\bigskip
\centerline{\bf Introduction}
\footnote{}{{\sifont Date}{\sfont : \vers}}

\bigskip\noindent
Let
$ Y $ be a complex hypersurface of degree
$ d $ in
$ X := \bP^{n} $ with
$ n\ge 2 $.
Let
$ f $ be a defining equation of
$ Y $, which is a homogeneous polynomial of degree
$ d $, and is assumed to be reduced.
Let
$ R $ be the Jacobian ring of
$ f $, and
$ R_{j} $ be its degree
$ j $ part.
If
$ Y $ is smooth, Griffiths' theorem [19] says that
$$
R_{qd+d-n-1} = H_{\rm prim}^{p,q}(Y,\bC),
$$
where
$ p+q = n-1 $ and the right-hand side is the
$ (p,q) $-part of the primitive cohomology, see also [18].
This follows from Griffiths' theorem saying that the Hodge filtration is
given by the pole order filtration in the smooth case, see [19].
The latter was generalized to the normal crossing case
by Deligne ([8], [9]), and the comparison between the Hodge
filtration and the pole order filtration was further studied in [10]
and [30] inspired by a letter of P.~Deligne which treated a certain
special case, see Remark 4.6 in [30].

Let
$ F $ be the Hodge filtration on
$ \cO_{X}(*Y) := \bigcup_{i>0}\cO_{X}(iY) $ as defined in [29].
It induces the Hodge filtration on the cohomology of the complement
$ U $ of
$ Y $ (see [9]) by using the induced filtration of the de Rham
complex, even if
$ Y $ is not a divisor with normal crossings.
Let
$ P $ denote the pole order filtration on
$ \cO_{X}(*Y) $ defined by
$ P_{i} = \cO_{X}((i+1)Y) $ if
$ i \ge 0 $ and
$ P_{i} = 0 $ otherwise.
Then we have
$ F_{i} \subset P_{i} $ in general, and
$ F_{i} = P_{i} $ for
$ i \le \alpha_{Y,y} - 1 $ on a neighborhood of
$ y \in \Sing\,Y $ if
$ -\alpha_{Y,y} $ is the maximal root of the Bernstein-Sato
polynomial divided by
$ s+1 $ for a local defining equation of
$ Y $ at
$ y $,
see [10], [30].
Set
$ \alpha_{Y} = \min_{y\in{\rm Sing}Y}\{\alpha_{Y,y}\} $.
We will denote also by
$ F $,
$ P $ the induced filtration on
$ \omega_{X}(*Y) = \omega_{X}\otimes_{\cO_{X}}\cO_{X}(*Y) $
(where
$ \omega_{X} = \Omega_{X}^{n} $).
Note that
$ F $,
$ P $ define the induced filtrations on the de Rham complex (see
(2.1.1)) and on the de Rham cohomology of
$ U $.

Let
$ \Omega^{\ssbull} $ denote the de Rham complex of global algebraic
differential forms on
$ \bA^{n+1} $.
Let
$ \ocH_{f} $ denote the free part of the algebraic Brieskorn module
$ \cH_{f} := \Omega^{n+1}/df\wdg d\Omega^{n-1} $, see [4].
It is known that
$ \ocH_{f} $ is a free graded
$ \bC[t] $-module of finite rank according to an algebraic version of [20]
(see also [1]), where the action of
$ t $ is given by the multiplication by
$ f $ and the degrees of
$ x_{i} $ and
$ dx_{i} $ are
$ 1 $.
(Here the
$ x_{i} $ are the coordinates of
$ \bA^{n+1} $.)
The action of the logarithmic Gauss-Manin connection
$ \nabla_{t\rd/\rd t} $ on the degree
$ k $ part
$ \ocH_{f,k} $ is the multiplication by
$ {k\over d}-1 $, and we have isomorphisms
$ t : \ocH_{f,k} \simto \ocH_{f,k+d} $ for
$ k \ge nd $, see Remark~(1.7)~(iii).
Let
$ \eta_{0} = {1\over d}\sum_{i=0}^{n}(-1)^{i}x_{i}dx_{0}\wdg
\cdots\widehat{dx_{i}}\cdots\wdg dx_{n} $.
For
$ q \in \bZ $, there are canonical morphisms
$$
H^{n}(U,\bC) \leftarrow
\Gamma(X,P_{q}(\omega_{X}(*Y))) \to \ocH_{f,(q+1)d},
\leqno(0.1)
$$
where the first morphism is defined
by using the de Rham cohomology and the second morphism sends
$ f^{-q-1}g\eta_{0} $ to
$ gdx_{0}\wdg\cdots\wdg dx_{n} $ for
$ g \in \bC[x]_{(q+1)d-n-1} $ and is surjective.
Note that
$ P_{q}(\omega_{X}(*Y)) $ is related to
$ P^{n-q} := P_{q-n} $ on
$ H^{n}(U,\bC) $ by definition of the filtration on the
de Rham complex, see (2.1.1).

\medskip\noindent
{\bf Theorem~1.} {\it
$ (i) $ For any
$ q \in \bZ $, the above morphisms in {\rm (0.1)} induce an isomorphism
$$
P^{n-q}H^{n}(U,\bC) = \ocH_{f,(q+1)d},
$$
such that the inclusion
$ P^{n-q+1}H^{n}(U,\bC) \hookrightarrow P^{n-q}H^{n}(U,\bC) $ is
identified with the morphism
$ t : \ocH_{f,qd} \hookrightarrow \ocH_{f,(q+1)d} $ and we have
$ P^{n-q}H^{n}(U,\bC) = H^{n}(U,\bC) $ for
$ q \ge n-1 $.

\medskip\noindent
$ (ii) $ The Hodge filtration
$ F^{n-q}H^{n}(U,\bC) $ is identified with the image of
$$
\Gamma(X,F_{q}(\omega_{X}(*Y))) \subset
\Gamma(X,P_{q}(\omega_{X}(*Y)))
$$
by the last morphism of {\rm (0.1)} for any
$ q \in \bZ $.
In particular,
}
$$
F^{n-q}H^{n}(U,\bC) = \ocH_{f,(q+1)d} \quad
\hbox{\it for}\,\,\, q \le \alpha_{Y} - 1.
$$

If
$ Y $ is a
$ \bQ $-homology manifold, then
$ F^{n-q}H^{n}(U,\bC) = F^{n-1-q}H_{\rm prim}^{n-1}(Y,\bC) $, and
$ H^{j}(Y,\bC) $ is naturally isomorphic to the intersection
cohomology
$ \IH^{j}(Y,\bC) $ (see [2], [15]).
If furthermore
$ Y $ is smooth, then
$ \alpha_{Y} = +\infty $, and
Theorem~1 is well known in the theory of hypersurface isolated
singularities, see for instance [33], [34], [35], [39].
The essential part of the proof of Theorem~1 is
showing the higher acyclicity of the Hodge filtration
on each component of the de Rham complex, see (2.2).
Note that it is not clear whether
$ F^{n-q} \ne P^{n-q} $ globally on
$ H^{n}(U,\bC) $ even if
$ F_{q} \ne P_{q} $ locally on
$ \cO_{X}(*Y) $, see (2.5).
Theorem~1 gives a refinement of the above mentioned theorem of
Griffiths by (1.6).
We can also generalize a description of the Kodaira-Spencer map
which is due to Griffiths in the smooth case, see (4.5).

If
$ Y $ is not smooth,
$ H^{j}(U,\bC) $ is identified, up to the non primitive part,
with the local cohomology
$ H^{j+1}_{Y}(X,\bC) $ or the homology
$ H_{2n-j-1}(Y,\bC) $.
For
$ q \ge [\alpha_{Y}] $, it is not easy to describe
$ F^{n-q}H^{n}(U,\bC) $ in general.
Consider first the case
$ q = 0 $.
Let
$ V $ denote the filtration on
$ \cO_{X} $ induced by the
$ V $-filtration of Kashiwara [22] and Malgrange [27] on
$ \cB_{f} = \cO_{X}\otimes_{\bC}\bC[\rd_{t}] $ along
$ Y $.
This coincides essentially with the filtration
by the multiplier ideals [24], see [5].
It induces a filtration
$ V $ on
$ \omega_{X}(Y) = \omega_{X}\otimes_{\cO}\cO_{X}(Y) $ such that
$ V^{\alpha}(\omega_{X}(Y)) = \omega_{X}(Y)\otimes_{\cO}
V^{\alpha+1}\cO_{X} $.
This is compatible with the filtration on
$ \omega_{X} $ defined by
$ V^{\alpha}\omega_{X} =\omega_{X}\otimes_{\cO}V^{\alpha}\cO_{X} $,
and gives a quotient filtration on
$ \omega_{Y} = \omega_{X}(Y)/\omega_{X} $.
Let
$ \tV^{>1}\cO_{X} $ be the inverse image of
$ V^{>0}\cB_{f} $ by the composition of the inclusion
$ \cO_{X}\to\cB_{f} $ with
$ \rd_{t} : \cB_{f}\to\cB_{f} $.
This gives
$ \tV^{>0}(\omega_{X}(Y)) := \tV^{>1}\omega_{X}\otimes
\cO_{X}(Y) $ and
$ \tV^{>0}\omega_{Y} := \tV^{>0}(\omega_{X}(Y))/\omega_{X} $
in the same way as above.

\medskip\noindent
{\bf Theorem~2.} {\it There are canonical isomorphisms
$$
\eqalign{
H^{0}(X,V^{0}(\omega_{X}(Y))) =
H^{0}(Y,V^{0}\omega_{Y}) &=
F^{n}H^{n}(U,\bC),\cr
H^{0}(X,\tV^{>0}(\omega_{X}(Y))) =
H^{0}(Y,\tV^{>0}\omega_{Y}) &=
F^{n-1}H^{n-1}(Y,\bC),\cr
}
$$
and the last term is canonically isomorphic to
$ F^{n-1}\IH^{n-1}(Y,\bC) $.
Furthermore,
$ \tV^{>1}\cO_{X} $ coincides with the adjoint ideal
{\rm ([14], [28], [38])}, and
$ \tV^{>0}\omega_{Y} $ is naturally isomorphic to the direct
image of the dualizing sheaf
$ \omega_{\tY} $ of a resolution of singularities
$ \tY $ of
$ Y $.
In particular,
$ F^{n}H^{n}(U,\bC) = F^{n-1}H^{n-1}(Y,\bC) $ in case
$ Y $ has only rational singularities.
}

\medskip
This follows from [30].
Note that
$ V^{0}\omega_{Y} = \omega_{Y} $ if and only if
$ \alpha_{Y} \ge 1 $, and
$ \tV^{>0}\omega_{Y} = \omega_{Y} $ if and only if
$ \alpha_{Y} > 1 $ (i.e.
$ Y $ has rational singularities).
See also [25] for the case of isolated singularities.

For
$ q \ge 1 $, it is not easy to describe
$ F_{q}(\cO_{X}(*Y)) $ unless we impose some strong condition on
the singularities of
$ Y $.
We now consider the case where
$ Y $ has only isolated singularities which are locally
semi-weighted-homogeneous.
This means that
$ Y $ is locally defined by a function
$ h = \sum_{\alpha\ge 1}h_{\alpha} $, where the
$ h_{\alpha} $ are weighted homogeneous polynomials of degree
$ \alpha $ for some appropriate local coordinates
$ y_{1},\dots,y_{n} $ with positive weights
$ w_{y,1},\dots,w_{y,n} $ around
$ y\in\Sing\,Y $, and
$ h_{1}^{-1}(0) $ (and hence
$ Y $) has an isolated singularity at
$ y $.
In this case, it is well known that
$ \alpha_{Y,y} = \sum_{i}w_{y,i} $, see (3.5) (i).
Let
$ \cO_{X,y}^{\ge \beta} $ be the ideal of
$ \cO_{X,y} $ generated by
$ \prod_{i}y_{i}^{\nu_{i}} $ with
$ \sum_{i}w_{y,i}\nu_{i} \ge \beta - \alpha_{Y,y} $.
Let
$ \cD_{X} $ be the sheaf of linear differential operators with the
filtration
$ F $ by the order of differential operators.
Put
$ k_{0} = [n-\alpha_{Y,y}] - 1 $.
Let
$ \cJ^{(q)} $ be the ideal sheaf of
$ \cO_{X} $ such that for
$ y \in \Sing\, Y $
$$
\cJ_{y}^{(q)}h^{-q-1} = \msum_{k=0}^{k_{0}}
F_{q-k}\cD_{X,y}(\cO_{X,y}^{\ge k+1}h^{-k-1}),
$$
and
$ \cJ^{(q)}|_{X\setminus{\rm Sing}\,Y} =
\cO_{X\setminus{\rm Sing}\,Y} $.
For example,
$ \cJ^{(0)} = \cO_{X}^{\ge 1} $ if
$ q = 0 $.
Let
$ J^{(q)} $ be the largest ideal of
$ \cO \,(= \cO_{\bA^{n+1}}) $ whose restriction to
$ \bA^{n+1}\setminus\{0\} $ coincides with the pull-back of
$ \cJ^{(q)} $ by the projection
$ \bA^{n+1}\setminus\{0\} \to \bP^{n} $.
Let
$ \ocH_{f}^{(q)} $ be the graded submodule of
$ \ocH_{f} $ defined by the image of
$ J^{(q)}dx_{0}\wdg\cdots\wdg dx_{n} $.
Then we have by [32] the following

\medskip\noindent
{\bf Theorem~3.} {\it With the above assumption,
there are canonical isomorphisms for
$ q \in \bN $
}
$$
\eqalign{
F_{q}(\cO_{X,y}(*Y)) &= \cJ_{y}^{(q)}h^{-q-1},\cr
F^{n-q}H^{n}(U,\bC) &= \ocH_{f,(q+1)d}^{(q)}.\cr}
$$

For example, if
$ q = 1 $,
$ w_{y,i} = 1/d_{y} $ for any
$ i $, and
$ \alpha_{Y,y} = n/d_{y} > 1 $ for any
$ y \in \Sing\,Y $, then
$ \cJ^{(1)} = \cO_{X}^{\ge 2} $.
Note that if
$ w_{y,i} = a_{y,i}^{-1} $ and the
$ a_{y,i} $ are integers which are mutually prime for any
$ y \in \Sing\,Y $, then
$ Y $ is a
$ \bQ $-homology manifold and the cohomology of
$ Y $ has good properties as remarked after Theorem~1.
We are informed that a formula similar to Theorem~3 is
obtained by L.~Wotzlaw if
$ Y $ has only ordinary double points as singularities and
$ n = 3 $.
See [36] for a completely different approach
to a similar problem.

In Sect.~1, we study meromorphic differential forms on projective
spaces for the proof of Theorem~1.
In Sect.~2, we prove the higher acyclicity of the Hodge filtration
on each component of the de Rham complex, and complete the proof of
Theorem~1.
In Sect.~3, we recall some facts from the theory of
$ V $-filtration, and prove Theorems 2 and 3.
In Sect.~4, we generalize a description of the Kodaira-Spencer
map.

\bigskip\bigskip
\centerline{\bf 1. Rational differential forms on projective spaces}

\bigskip\noindent
{\bf 1.1.}
Let
$ f $ be a reduced homogeneous polynomial of degree
$ d $, and put
$ X = \bP^{n} $,
$ Y = f^{-1}(0) \subset X $, and
$ U = X \setminus Y $.
Let
$ x_{0}, \dots, x_{n} $ be the coordinates of
$ \bA^{n+1} $, and set
$$
\xi = {1\over d}\msum_{i} x_{i}\rd_{i},
$$
so that
$ \xi f = f $, where
$ \rd_{i} = \rd/\rd x_{i} $.
Let
$ \iota_{\xi} $ and
$ L_{\xi} $ denote respectively the interior product and the
Lie derivation.
Then
$$
\iota_{\xi}(\dfoverf) = 1,\quad
\iota_{\xi}\scirc\iota_{\xi} = 0,\quad
\iota_{\xi}\scirc d + d\scirc\iota_{\xi} = L_{\xi}.
\leqno(1.1.1)
$$

Let
$ \Omega^{\ssbull} $ denote the de Rham complex of global
differential forms on
$ \bA^{n+1} $.
This is isomorphic to the Koszul complex for the action of the
vector fields
$ \rd_{i} $ on the polynomial ring
$ \bC[x_{0}, \dots, x_{n}] $.
Let
$$
\cH_{f} = \Omega^{n+1}/df\wdg d\Omega^{n-1},
$$
as in the introduction.
It has a structure a graded module
$ \cH_{f} = \bigoplus_{k\in\bZ}\cH_{f,k} $ where the degree of
$ x_{i} $ and
$ dx_{i} $ are
$ 1 $.
Let
$ \Omega^{\ssbull}[f^{-1}]_{k} $ denote the degree
$ k $ part of the localization
$ \Omega^{\ssbull}[f^{-1}] $ of
$ \Omega^{\ssbull} $.
Then the restriction of
$ L_{\xi} $ to
$ \cH_{f,k} $ and
$ \Omega^{\ssbull}[f^{-1}]_{k} $ is the multiplication by
$ k/d $.
Let
$ \Omega^{\ssbull}[f^{-1}]_{0}^{(\xi)} $ be the subcomplex of
$ \Omega^{\ssbull}[f^{-1}]_{0} $ defined by
$$
\Omega^{j}[f^{-1}]_{0}^{(\xi)} = \Im(\iota_{\xi} :
\Omega^{j+1}[f^{-1}]_{0} \to \Omega^{j}[f^{-1}]_{0}).
$$

\medskip\noindent
{\bf 1.2.~Lemma.} {\it
There is a canonical isomorphism of complexes}
$$
\Gamma(U,\Omega_{U}^{\ssbull}) =
\Omega^{\ssbull}[f^{-1}]_{0}^{(\xi)}.
\leqno(1.2.1)
$$

\medskip\noindent
{\it Proof.}
This follows by using the blow-up of
$ \bP^{n+1} $ along the origin of
$ \bA^{n+1}\,(\subset \bP^{n+1}) $ together with the pull-back
by the projection to
$ \bP^{n} = \bP^{n+1}\setminus\bA^{n+1} $.
See also Proposition 6.1.16 in [12].

\medskip\noindent
{\bf 1.3.~Lemma.} {\it
The canonical morphism
$$
\dfoverf\wdg : \Omega^{\ssbull}[f^{-1}]_{0}^{(\xi)} \to
\dfoverf\wdg\Omega^{\ssbull}[f^{-1}]_{0}^{(\xi)}
\leqno(1.3.1)
$$
is an isomorphism of complexes, and its inverse is given by
}
$$
\iota_{\xi}:\dfoverf\wdg\Omega^{\ssbull}[f^{-1}]_{0}^{(\xi)}\simto
\Omega^{\ssbull}[f^{-1}]_{0}^{(\xi)}.
\leqno(1.3.2)
$$

\medskip\noindent
{\it Proof.}
By (1.1.1), we have
$ \iota_{\xi}\scirc(\dfoverf\wdg) = id $ on
$ \Omega^{\ssbull}[f^{-1}]_{0}^{(\xi)} $,
and (1.3.1) is surjective by definition.

\medskip\noindent
{\bf 1.4.~Lemma.} {\it
The canonical morphism
$$
H^{j}(\dfoverf\wdg\Omega^{\ssbull}[f^{-1}]_{0}^{(\xi)})\to
H^{j}(\dfoverf\wdg\Omega^{\ssbull}[f^{-1}]_{0})
\leqno(1.4.1)
$$
is injective for any
$ j $, and is bijective for
$ j = n $.
}

\medskip\noindent
{\it Proof.}
Take
$ f^{-i}\eta\in\Omega^{j-2}[f^{-1}]_{0} $.
Here we may assume
$ i > 0 $, replacing
$ \eta $ with
$ f^{k}\eta $ and
$ i $ with
$ i + k $.
Since the degree of
$ \eta $ is
$ di $, we have
$ i\eta = L_{\xi}\eta = \iota_{\xi}(d\eta)+d(\iota_{\xi}\eta) $,
and
$$
\dfoverf\wdg d(if^{-i}\eta) =
\dfoverf\wdg d(\iota_{\xi}(f^{-i}d\eta)).
$$
So the injectivity follows.
The surjectivity follows from
$$
(\dfoverf\wdg)\scirc\iota_{\xi} = id\quad\hbox{on}\,\,\,
\Omega^{n+1}[f^{-1}].
\leqno(1.4.2)
$$

\medskip\noindent
{\bf 1.5.~Proposition.} {\it
There are canonical isomorphisms
$$
H^{n}(\Gamma(U,\Omega_{U}^{\ssbull})) =
H^{n}(\dfoverf\wdg\Omega^{\ssbull}[f^{-1}]_{0}^{(\xi)}) =
\mathop{\smash{\mathop{\hbox{\rm lim}}\limits_{\raise.4ex
\hbox{${\scriptstyle\longrightarrow}$}}}}
\limits_{\smash{\raise-0.7ex\hbox{$\scriptstyle k$}}}
\cH_{f,kd},
$$
where the first isomorphism is induced by {\rm (1.2.1)} and
{\rm (1.3.1)}, and the limit is taken over
$ k \in \bN $ with transition morphisms
$ \cH_{f,kd} \to \cH_{f,(k+1)d} $ given by the multiplication by
$ f $.
}

\medskip\noindent
{\it Proof.}
The first isomorphism follows from Lemmas (1.2) and (1.3).
The second isomorphism is defined by assigning
$ \eta $ to
$ \eta/f^{k} \in \Omega[f^{-1}]_{0} $, because (1.4.2) implies
$$
\dfoverf\wdg\Omega^{n}[f^{-1}]_{0}^{(\xi)} =
\Omega^{n+1}[f^{-1}]_{0}.
\leqno(1.5.1)
$$
Since
$ \eta/f^{k} =  f\eta/f^{k+1} $, the isomorphism then follows from
the definition of
$ \cH_{f} $ together with Lemma (1.4).

\medskip\noindent
{\bf 1.6.~Proposition.} {\it
We have a canonical isomorphism
$$
\cH_{f}/f\cH_{f} = \Omega^{n+1}/df\wdg\Omega^{n} \,
(= \cO/(\rd f) =: R),
\leqno(1.6.1)
$$
where
$ (\rd f) $ is the Jacobian ideal generated by the
$ f_{i} := \rd f/\rd x_{i} $ in
$ \cO = \bC[x_{0},\dots,x_{n}] $.
}

\medskip\noindent
{\it Proof.}
The assertion is equivalent to
$$
df\wdg d\Omega^{n-1} + f\Omega^{n+1} = df\wdg\Omega^{n}\,
(=(\rd f)\Omega^{n+1}).
$$
Since
$ f = \sum_{i}{1\over d}x_{i}\rd_{i}f/\rd x_{i} $, the inclusion
$ \subset $ is clear.
For the converse we use the Gauss-Manin connection
$ \nabla_{\rd_{t}} $ on
$ \cH_{f} $.
(We do not know a simple proof without using essentially the
Gauss-Manin connection even in the isolated singularity case.)
We have by definition
$$
\nabla_{\rd_{t}}\omega = d\eta \quad\hbox{with}\quad
df\wdg\eta = \omega,
$$
where
$ \rd_{t} = \rd/\rd t $.
Note that the inverse of the Gauss-Manin connection
$ \nabla_{\rd_{t}}^{-1} $ is well-defined as a
$ \bC $-linear endomorphism of
$ \cH_{f} $ by the de Rham lemma.
It is well known that
$$
\nabla_{\rd_{t}}(f\omega) = (k/d)\omega\quad\hbox{for}\quad
\omega \in \cH_{f,k}.
\leqno(1.6.2)
$$
Indeed, this follows from
$ d(\iota_{\xi}\omega) = L_{\xi}\omega $ (see (1.1.1)) together
wtih
$ \dfoverf\wedge(\iota_{\xi}\omega) = \omega $ by setting
$ \eta = \iota_{\xi}\omega $.

We have to show that
$ \omega \in f\cH_{f} $ if
$ \omega $ is represented by an element of
$ df\wdg\Omega^{n} $.
Here we may assume
$ \omega \in \cH_{f,k} $ (because
$ f $ is homogeneous), and we have
$ k > d $ by the definition of the grading on
$ \Omega^{n+1} $ (because
$ \deg f_{i} = d-1 $ and
$ \deg dx_{0}\wdg\cdots\wdg dx_{n} = n+1 > 1 $ ).
So the assertion follows from (1.6.2) together with
$ \nabla_{\rd_{t}}t = t\nabla_{\rd_{t}} + id $.

\medskip\noindent
{\bf 1.7.~Remarks.}
(i) Let
$ {}_{\rm tor}\cH_{f} $ denote the torsion part of
$ \cH_{f} $ as a
$ \bC[t] $-module so that we have a short exact sequence
$$
0 \to {}_{\rm tor}\cH_{f} \to \cH_{f} \to \ocH_{f} \to 0.
$$
It is well known after [20] that
$ \ocH_{f} $ is finitely generated over
$ \bC[t] $ and
$ {}_{\rm tor}\cH_{f} $ is killed by a sufficiently high power
of
$ f $.
(Indeed, this is easily proved by using a resolution of
singularities, see e.g. [1], 2.3.)
Then, by Proposition (1.6) together with the snake lemma
applied to the multiplication by
$ f $ on the above exact sequence, we get
$$
{}_{\rm tor}\cH_{f,k}/f\,{}_{\rm tor}\cH_{f,k-d} =
\cH_{f,k}/f\cH_{f,k-d} = R_{k-n-1}\quad\hbox{for}\quad k\gg 0.
$$
In particular,
$ {}_{\rm tor}\cH_{f} $ is not finitely generated over
$ \bC[t] $ unless
$ Y $ is smooth.
In the case where
$ Y $ has only isolated singularities, then the dimension of
$ R_{k-n-1} $ is closely related to the Tjurina numbers,
see [7].

\medskip
(ii) In [1] the Brieskorn module is defined by the cohomology
$ H^{j}\cA_{f}^{\ssbull} $ of
the complex
$ (\cA_{f}^{\ssbull},d) $ where
$ \cA_{f}^{j} = \Ker(df\wdg:\Omega^{j}\to\Omega^{j+1}) $.
For
$ j = n + 1 $, we have a canonical surjection
$$
H^{n+1}\cA_{f}^{\ssbull}\to\cH_{f}\,
(=\Omega^{n+1}/df\wdg d\Omega^{n-1}),
$$
and its kernel is
$ f $-torsion by the acyclicity of
$ (\Omega^{\ssbull}[f^{-1}],df\wdg) $.
So the definition in [1] is compatible with the one in this paper
as long as we take the free part.
(However, it is not clear whether Proposition (1.6) holds with
$ \cH_{f} $ replaced by
$ H^{n+1}\cA_{f}^{\ssbull} $.)

\medskip
(iii) For
$ k \ge nd $, we have isomorphisms
$ t : \ocH_{f,k} \simto \ocH_{f,k+d} $.
This follows for example from Remark (3.5) (ii) below
because (3.5.1) implies the bijectivity of
$ t : \Gr_{V}^{\alpha}\ocH_{f}\to\Gr_{V}^{\alpha+1}\ocH_{f} $
for
$ \alpha > \beta_{f} - 1 $ and we have
$ \beta_{f} - 1 < n $.
We can also show the isomorphism by proving Theorem~1 with
$ q \ge n $ in the first assertion replaced by
$ q \gg n $, and then using the fact that the algebraic de Rham
cohomology of
$ U $ is generated by meromorphic differential forms having
poles of order
$ \le n $ along
$ Y $, which follows from the relation between the Hodge and
the pole order filtrations, see [10].

\medskip\noindent
{\bf 1.8.~Milnor cohomology.}
With the notation of (1.1), let
$ \cS^{i} = \cO_{X}(-i) $ and
$ E = \cSpec_{X}(\mopls_{i\in\bN}\cS^{i}) $.
Then
$ E $ is a line bundle over
$ X $, and the sheaf of local sections
$ \cO_{X}(E) $ is naturally identified with
$ \cO_{X}(1) $.
Let
$ f $ be as in (1.1).
Since it is identified with a section of
$ \cO_{X}(d) $, it induces a shifted graded morphism
$ \rho_{f} : \cS^{i}\to \cS^{i-d} $ for
$ i \ge d $, which is compatible with the action of
$ u \in \cS^{j} $ defined by the multiplication
$ \cS^{j}\otimes\cS^{i}\to\cS^{i+j} $.
So it determines an ideal
$ \cJ $ of
$ \cS^{\ssbull} := \mopls_{i\in\bN}\cS^{i} $ which is locally
generated by
$ u - \rho_{f}(u) $ for local nonzero sections
$ u $ of
$ \cO_{X}(-d) $.
Let
$ \ocS = \cS^{\ssbull}/\cJ $, and
$ Z = \cSpec_{X}\ocS $ with the canonical projection
$ \pi : Z \to X $.
Then
$ Z $ is a divisor on
$ E $, and we have as
$ \cO_{X} $-modules
$$
\pi_{*}\cO_{Z} = \ocS = \mopls_{0\le i<d}\cS^{i}.
\leqno(1.8.1)
$$

Let
$ F = f^{-1}(1) \subset \bA^{n+1} $, and
$ \oF $ be the closure of
$ F $ in the compactification
$ \bP^{n+1} $ of
$ \bA^{n+1} $.
Then
$ Z $ can be identified with
$ \oF $.
Indeed, we have a projection
$ \pi' : \oF \to X $ induced by the projection
$ \bP^{n+1}\setminus \{0\}\to X = \bP^{n} $ (where
$ 0 \in \bA^{n+1} \subset \bP^{n+1} $), and
$ F = \pi'{}^{-1}(U) $.
Furthermore,
$ \bA^{n+1}\setminus \{0\} $ is identified with the total space
of the
$ \bC^{*} $-bundle associated to
$ \cO_{X}(-1) $, and
$ \bP^{n+1}\setminus \{0\} $ is identified with the line bundle
corresponding to
$ \cO_{X}(1) $ by exchanging the zero section and the
$ \infty $-section (using the involution of
$ \bC^{*} $ sending
$ \lambda $ to
$ \lambda^{-1} $).
So we get a canonical isomorphism
$$
(\oF,F) = (Z,\pi^{-1}(U))\quad\hbox{over}\,\,\, U.
\leqno(1.8.2)
$$

Let
$ G = \bZ/d\bZ $ be the covering transformation group of
$ F\to U $ and
$ Z\to X $ such that
$ F/G = U $ and
$ Z/G = X $.
It has a generator
$ g $ which acts on
$ F $ by
$$
g : (x_{0},\dots,x_{n}) \mapsto (\zeta x_{0},\dots,\zeta x_{n}),
$$
where
$ \zeta = \exp(2\pi i/d) $.
This coincides with the geometric Milnor monodromy.
Furthermore, the subsheaf
$ \cS^{i} $ of
$ \pi_{*}\cO_{Z} $ for
$ 0 \le i < d $ in
(1.8.1) coincides with the eigenspace
corresponding to the eigenvalue
$ \zeta^{i} $ for the action of
$ g^{*} $ because the above involution of
$ \bC^{*} $ is used in the isomorphism (1.8.2).
In particular,
$ \cS^{i}|_{U} $ is stable by the Gauss-Manin connection.

Using these we can generalize the first assertion of Theorem~1
to the Milnor cohomology.
Indeed, Proposition (1.5) is generalized to the following
assertion for
$ 0 \le i < d $:
$$
H^{n}(\Gamma(U,(\cS^{i}|_{U})\otimes\Omega_{U}^{\ssbull})) =
H^{n}(\dfoverf\wdg\Omega^{\ssbull}[f^{-1}]_{-i}^{(\xi)}) =
\mathop{\smash{\mathop{\hbox{\rm lim}}\limits_{\raise.4ex
\hbox{${\scriptstyle\longrightarrow}$}}}}
\limits_{\smash{\raise-0.7ex\hbox{$\scriptstyle k$}}}
\cH_{f,kd-i}.
\leqno(1.8.3)
$$
For
$ 0 \le i < d $, we define the pole order filtration on
$ \cS_{i}\otimes\cO_{X}(*Y) $ by
$$
P_{j}(\cS_{i}\otimes\cO_{X}(*Y)) =
\cS_{i}\otimes\cO_{X}((j+1)Y)\quad\hbox{if}\,\,\, j \ge 0,
$$
and it is zero otherwise.
We have the Hodge filtration
$ F $ on
$ \cS_{i}\otimes\cO_{X}(*Y) $ because (1.8.1) is the
decomposition by the action of
$ g $.
Furthermore,
$ F_{j}\subset P_{j} $, and they coincide at the smooth point of
$ D $.

For
$ q \ge n+1 $, there are canonical isomorphisms
$ H^{n}(F,\bC)_{\be(i/d)} = \ocH_{f,qd-i} $ compatible with
$ t : \ocH_{f,qd-i} \simto \ocH_{f,(q+1)d-i} $, where
$ H^{n}(F,\bC)_{\be(i/d)} $ is the
$ \be(i/d) $-eigenspace for the action of the monodromy, where
$ \be(i/d) = \exp(2\pi\sqrt{-1}i/d) $.
For
$ q \in \bZ $, we have the surjective morphisms
$$
\Gamma(X,P_{q}(\cS_{i}\otimes\omega_{X}(*Y))) \to
\ocH_{f,(q+1)d-i},
\leqno(1.8.4)
$$
and, for
$ q \le n $,
$ P^{n-q}H^{n}(F,\bC)_{\be(i/d)} $ coincides with the image of
$$
t^{n-q}: \ocH_{f,(q+1)d-i} \hookrightarrow
\ocH_{f,d(n+1)-i} = H^{n}(F,\bC)_{\be(i/d)},
$$
However, the remaining assertion of Theorem~1 cannot be
generalized, because we have to use the shifted
$ b $-function
$ b_{f}(s+i/d) $ corresponding to
$ f^{i/d}f^{s} = f^{s+i/d} $, and consider the minimal root of
$ b_{f}(-s+i/d) $ instead of
$ b_{f}(-s+i/d)/(-s+i/d+1) $.

\medskip\noindent
{\bf 1.9.~Remark.}
Let
$ \omega_{0} = dx_{0} \wdg \cdots \wdg dx_{n} $, and
consider the following Brian\c con-Skoda type property:
$$
\hbox{There is a positive integer
$ k $ such that}\,\,\,
f^{k} \omega_{0} \in df \wdg d\Omega^{n-1}.
$$
This property is clearly equivalent to
$ [\omega_{0}] = 0 $ in
$ \ocH_{f,n+1} $.
By the above discussion, we get the following corollary:

\medskip
Assume
$ e(-(n+1)/d) $ is not an eigenvalue of the monodromy acting on
$ H^{n}(F,\bC) $.
Then
$ f $ satisfies the above Brian\c con-Skoda type property.

\medskip
The above hypothesis is satisfied, for example, if
$ f = x^{3}+y^{2}z $, the equation of a cuspidal cubic plane
curve, or
$ f = x^{2}z+y^{3}+xyt $, the equation of a cubic surface such that
$ H^{3}(F)=0 $, see for details [11], Example 4.3.
Higher dimensional examples can be constructed easily using the
hypersurfaces introduced in [12], p. 148, Proposition (2.24).

\bigskip\bigskip
\centerline{\bf 2. Hodge and pole order filtrations}

\bigskip\noindent
{\bf 2.1} With the notation of (1.1), let
$ \cO_{X}(*Y) $ be the localization of
$ \cO_{X} $ by local defining equations of
$ Y $.
Let
$ P $ and
$ F $ be respectively the pole order filtration and the Hodge
filtration on the left
$ \cD_{X} $-module
$ \cO_{X}(*Y) $ (see [10], [29]) as in the introduction (e.g.
$ P_{i} = \cO_{X}((i+1)Y) $ if
$ i \ge 0 $ and
$ P_{i} = 0 $ otherwise).
They induce the filtrations
$ F $ and
$ P $ on the de Rham complex
$ \DR_{X}(\cO_{X}(*Y)) = \Omega_{X}^{\ssbull}(*Y) $ by
$$
F_{i}(\Omega_{X}^{j}(*Y)) =
\Omega_{X}^{j}\otimes F_{i+j}(\cO_{X}(*Y)),
\leqno(2.1.1)
$$
and similarly for
$ P $.
Note that the filtrations
$ F $,
$ P $ on
$ \Omega_{X}^{n} $ are different from those on
$ \omega_{X} $ defined in the introduction (i.e. there is a
shift by
$ n $).

If
$ i > 0 $, it follows from Bott's vanishing theorem that
$$
H^{k}(X,\Omega_{X}^{j}(iY)) = 0\quad\hbox{for}\quad k > 0,
\leqno(2.1.2)
$$
This implies the
$ \Gamma $-acyclicity of the components of
$ P_{i}\DR_{X}(\cO_{X}(*Y)) $, i.e.
$$
H^{k}(X,P_{i}(\Omega_{X}^{j}(*Y))) = 0\quad\hbox{for}\quad k>0.
\leqno(2.1.3)
$$

\medskip\noindent
{\bf 2.2.~Proposition.} {\it
$ H^{k}(X,F_{i}(\Omega_{X}^{j}(*Y))\otimes_{\cO_{X}}\cO_{X}(r))
= 0 $ for
$ k > 0 $,
$ r \ge 0 $ and
$ i,j \in \bZ $.
}

\medskip\noindent
{\it Proof.}
We proceed by increasing induction on
$ n = \dim X \ge 0 $ and increasing induction on
$ i \ge -n $.
The assertion is trivial if
$ n = \dim X = 0 $.
Since
$ F_{-1}\cO_{X}(*Y) = 0 $, we have
$$
F_{-n}(\Omega_{X}^{\ssbull}(*Y)) =
\Omega_{X}^{n}\otimes F_{0}(\cO_{X}(*Y))[-n].
$$
Combining this with the vanishing of
$ H^{j}(U,\bC) $ for
$ j > n $,
the assertion for
$ i = -n $ and
$ r = 0 $ follows from the strictness of the
Hodge filtration on the direct image of
$ (\cO_{X}(*Y),F) $ by
$ X \to pt $.
Indeed, the latter means the injectivity of
$$
H^{j}(X,F_{p}(\Omega_{X}^{\ssbull}(*Y))) \to
H^{j}(X,\Omega_{X}^{\ssbull}(*Y)) = H^{j}(U,\bC),
\leqno(2.2.1)
$$
(even if
$ Y $ is not a divisor with normal crossings),
because the direct image by
$ X \to pt $ is defined by using the de Rham
complex.

Let
$ X_{0} $ be a general hyperplane of
$ X = \bP^{n} $, and set
$ Y_{0} = Y \cap X_{0} $.
Then
$ X_{0} $ is non characteristic to the
$ \cD_{X} $-module
$ \cO_{X}(*Y) $, and the vanishing cycle
$ \cD $-module
$ \varphi_{g}\cO_{X}(*Y) \,(= \mopls_{0\le\alpha<1}\Gr_{V}^{\alpha}
\cO_{X}(*Y)) $ vanishes, where
$ g $ is a local equation of
$ X_{0} \subset X $.
This implies that the filtration
$ V $ on
$ \cO_{X}(*Y) $ along
$ X_{0} $ is given by the
$ g $-adic filtration (because
$ \Gr_{V}^{\alpha}\cO_{X}(*Y) = 0 $ unless
$ \alpha $ is a positive integer).
So the restriction of the filtered (left)
$ \cD_{X} $-module
$ (\cO_{X}(*Y) ,F) $ by the inclusion
$ i_{0} : X_{0} \to X $ is given by the tensor product with
$ \cO_{X_{0}} $ over
$ \cO_{X} $ (because this gives
$ \Gr_{V}^{1} $).
Then, by the uniqueness of the direct image of a mixed Hodge
module by the affine open morphism
$ X_{0}\setminus Y_{0} \to X_{0} $ (see [29], 2.11), we get
$$
F_{i}(\cO_{X_{0}}(*Y_{0})) =
F_{i}(\cO_{X}(*Y))\otimes_{\cO}\cO_{X_{0}}.
$$
Since
$ \Omega_{\bP^{n}}^{n} = \cO_{\bP^{n}}(-n-1) $,
we get furthermore a short exact sequence
$$
0 \to F_{i}(\Omega_{X}^{n}(*Y))(r) \to F_{i}(\Omega_{X}^{n}(*Y))(r+1)
\to F_{i}(\Omega_{X_{0}}^{n-1}(*Y_{0}))(r) \to 0,
$$
where
$ M(r) $ denotes
$ M \otimes_{\cO}\cO_{X}(r) $ for an
$ \cO_{X} $-module
$ M $.
Using the inductive hypothesis for
$ X_{0} $, this implies
$$
\eqalign{
H^{k}(X,F_{i}(\Omega_{X}^{n}(*Y))(r+1))
& = 0\quad\hbox{for any}\,\,\,k > 0, \cr
\hbox{if}\quad
H^{k}(X,F_{i}(\Omega_{X}^{n}(*Y))(r))
& = 0\quad\hbox{for any}\,\,\,k > 0. \cr
}
\leqno(2.2.2)
$$

For
$ n \ge 1 $, there is a short exact sequence
$$
0\to\cO_{X}\to\mathop{\mopls}\limits^{n+1}\cO_{X}(1)\to
\Theta_{X}\to 0,
$$
where
$ \Theta_{X} $ is the sheaf of vector fields.
This induces a short exact sequence
$$
0\to\bwdg^{i-1}\Theta_{X}\to\bwdg^{i}
(\mathop{\mopls}\limits^{n+1}\cO_{X}(1))
\to\bwdg^{i}\Theta_{X}\to 0.
$$
Since
$ \Omega_{X}^{j} = \Omega_{X}^{n}\otimes\bwdg^{n-j}\Theta_{X} $,
it implies by decreasing induction on
$ j < n $
$$
\eqalign{
H^{k}(X,F_{i+n-j}(\Omega_{X}^{j}(*Y))(r))
& = 0\quad\hbox{for any}\,\,\,k > 0,\, r \ge 0, \cr
\hbox{if}\quad
H^{k}(X,F_{i}(\Omega_{X}^{n}(*Y))(r))
& = 0\quad\hbox{for any}\,\,\,k > 0,\, r \ge 0. \cr
}
\leqno(2.2.3)
$$
Here the hypothesis is reduced to the case
$ r = 0 $ by (2.2.2).

Consider now the spectral sequence
$$
E_{1}^{p,q} = H^{q}(X,F_{i}(\Omega_{X}^{p}(*Y))) \Rightarrow
H^{p+q}(X,F_{i}(\Omega_{X}^{\ssbull}(*Y))).
\leqno(2.2.4)
$$
By inductive hypothesis for
$ i $ together with (2.2.2) and (2.2.3), we have
$ E_{1}^{p,q} = 0 $ for
$ p < n $,
$ q > 0 $.
On the other hand,
$ H^{k}(X,F_{i}(\Omega_{X}^{\ssbull}(*Y))) = 0 $ for
$ k > n $ by the strictness of the Hodge filtration.
So
$ E_{1}^{n,q} = H^{q}(X,F_{i}(\Omega_{X}^{n}(*Y))) = 0 $ for
$ q > 0 $, and we can proceed by increasing induction on
$ i $ using (2.2.2) and (2.2.3) (with
$ n $ fixed).
This completes the proof of Proposition (2.2).

\medskip\noindent
{\bf 2.3.~Corollary.} {\it We have canonical isomorphisms for
$ q \in \bZ $
}
$$
F^{n-q}H^{n}(U,\bC) = \Im(\Gamma(X,F_{q}(\Omega_{X}^{n}(*Y))\to
H^{n}(X,\Omega_{X}^{\ssbull}(*Y))).
$$

\medskip\noindent
{\it Proof.}
By the definition of the direct image of mixed Hodge modules [29],
we have a canonical isomorphism
$$
F^{n-q}H^{n}(U,\bC) = \Im(H^{n}(X,F_{q-n}(\Omega_{X}^{\ssbull}(*Y)))
\to H^{n}(X,\Omega_{X}^{\ssbull}(*Y))).
$$
So the assertion follows from Proposition (2.2).

\medskip\noindent
{\bf 2.4.~Proof of Theorem~1.}
The assertion (i) follows from (1.5).
We have
$ F_{q}\cO_{X} = P_{q}\cO_{X} $ for
$ q \le \alpha_{Y,y} - 1 $ on a neighborhood of
$ y \in \Sing\,Y $ as in the introduction.
Since
$ \ocH_{f} $ is the quotient of
$ \cH_{f} $ by the
$ f $-torsion part, we see that
$ \ocH_{f,j} $ is isomorphic to the image of
$ \cH_{f,j} $ in the inductive limit of
$ \cH_{f,j+kd} $ over
$ k \in \bN $ where the transition morphisms are given by
the multiplication by
$ f $.
So the assertion (ii) follows from (2.3).

\medskip\noindent
{\bf 2.5.~Remark.}
If
$ \alpha_{Y,y} < 1 $, we have locally
$ F_{0} \ne P_{0} $ on
$ \cO_{X,y}(*Y) $.
However, it is not clear whether
$ F^{n} \ne P^{n} $ on
$ H^{n}(U,\bC) $.
For example, in the case
$ n = 2 $, we have
$ F^{2} = P^{2} $ on
$ H^{2}(U,\bC) $ if
$ d \le 3 $.
Indeed, if
$ Y $ has only ordinary double points and
$ n = 2 $, then
$ \alpha_{Y} = 1 $ and hence
$ F_{0} = P_{0} $.
In case
$ Y $ is a rational cubic curve with a cusp (i.e.
$ f = x^{3} + y^{2}z $), we have
$ \alpha_{Y} = 5/6 $ but
$ H^{2}(U,\bC) = 0 $ because
$ Y $ is a
$ \bQ $-rational manifold.
However, for $ d = 4 $, we have
$ F^{2} \ne P^{2} $ on
$ H^{2}(U,\bC) $ if
$ Y $ is the union of a smooth cubic curve and a line
which intersect only at one point,
see [12], p.~186, Remark 1.33.
We have also
$ F^{2} \ne P^{2} $ if
$ d = 4 $ and
$ Y $ has two cusps
$ O, O' $ so that its normalization is an elliptic curve
$ E $, e.g. if
$ f = x^{2}y^{2} + xz^{3} + yz^{3} $.
Indeed, let
$ g_{1}, g_{2} $ be linear functions of
$ x, y, z $ such that
$ g_{1}^{-1}(0) $ passes through both
$ O, O' $ and
$ g_{2}^{-1}(0) $ passes through
$ O $ but not
$ O' $.
Let
$ \omega_{i} $ be the differential
$ 2 $-form on
$ U $ corresponding to
$ f^{-1}g_{i}dx\wedge dy\wedge dz $ by Theorem~1, and
$ \eta_{i} $ be the differential
$ 1 $-form on
$ Y \setminus \{O,O'\} $ obtained by taking the residue of
$ \omega_{i} $.
Let
$ E $ be the normalization of
$ Y $.
Then
$ \eta_{1} $ is extended to a nowhere vanishing
$ 1 $-form on
$ E $ by Theorem~3, and hence
$ \eta_{2} $ has only a pole of order
$ 2 $ at the point
$ P' $ corresponding to
$ O' $ because the pull-back of
$ g_{2}/g_{1} $ to
$ E $ has only such a pole.
Since
$ E':=E\setminus\{P'\} $ is affine, the cohomology of
$ E' $ is calculated by the complex of algebraic differential
forms on
$ E' $, and the cohomology class of
$ \eta_{2} $ modulo
$ \bC\eta_{1} $ does not vanish, i.e. it is not an exact form
modulo
$ \bC\eta_{1} $, because there is no rational function
on an elliptic curve having only one pole of order
$ 1 $.

Note however that we have
$ F = P $ on
$ H^{2}(U,\bC) $ if any singular point of
$ Y $ is a cusp or an ordinary double point and if
the irreducible components of
$ Y $ are rational curves.
Indeed, we have
$ F^{2}H^{2}(U,\bC) = H^{2}(U,\bC) $ using the weight spectral
sequence because the normalization of each irreducible component
of
$ Y $ is a smooth rational curve.

If
$ \dim Y = n - 1 = 2 $ and
$ Y $ has only rational singularities, then we have locally
$ F_{1} \ne P_{1} $ on
$ \cO_{X}(*Y) $ by [32], but it is not easy to show that
$ F^{2} \ne P^{2} $ globally on
$ H^{3}(U,\bC) $, or equivalently,
$ F^{1} \ne P^{1} $ on
$ H^{2}(Y,\bC) $.
In this case,
$ Y $ is a
$ \bQ $-homology manifold and if
$ Y' $ is a smooth hypersurface of degree
$ d $ in
$ \bP^{3} $, then we have
$ h^{p,p-2}(Y') = h^{p,p-2}(Y) $ for
$ p \ne 1 $ and
$ h^{1,1}(Y') - h^{1,1}(Y) $ is the sum of the Milnor numbers.
More generally, if
$ Y $ is a
$ \bQ $-homology manifold, then
$ h^{p,q}(Y) $ may be expressed in terms of the degree
$ d $ and the list of singularities on
$ Y $.
For details, see [13].

\bigskip\bigskip
\centerline{\bf 3. $ V $-Filtration}

\bigskip\noindent
{\bf 3.1.}
Let
$ X $ be a smooth complex algebraic variety, and
$ h $ be an algebraic function on
$ X $.
Set
$ Y = h^{-1}(0) $.
Let
$ i_{h} : X \to X\times\bA^{1} $ be the graph embedding.
We denote by
$ \cB_{h} $ the direct image of
$ \cO_{X} $ by
$ i_{h} $ as a
$ \cD $-module, see [30] for more details.
Let
$ \cB_{h}[t^{-1}] $ be the localization of
$ \cB_{h} $ by
$ t $.
This is identified with the direct image of
$ \cO_{X}(*Y) $ by
$ i_{h} $ as a
$ \cD $-module.

We have canonical isomorphisms
$$
\cB_{h} = \cO_{X}\otimes_{\bC}\bC[\rd_{t}],\quad
\cB_{h}[t^{-1}] = \cO_{X}(*Y)\otimes_{\bC}\bC[\rd_{t}],
$$
where
$ t $ is the coordinate of
$ \bA^{1} $ and
$ \rd_{t} = \rd/\rd t $.
These isomorphisms are compatible with the action of
$ \cO_{X} $ and
$ \rd_{t} $.
The actions of
$ t $ and of a vector field
$ \xi $ on
$ X $ are given by
$$
\eqalign{
\xi(g\otimes\rd^{i})
&= (\xi g)\otimes\rd^{i} - (\xi h)g\otimes\rd^{i+1},\cr
t(g\otimes\rd^{i})
&= hg\otimes\rd^{i} - ig\otimes\rd^{i-1}.\cr
}
\leqno(3.1.1)
$$

We have the Hodge filtration
$ F $ on
$ \cB_{h} $,
$ \cB_{h}[t^{-1}] $ by
$$
F_{p}(\cB_{h}[t^{-1}]) = \msum_{i\ge 0}
F_{p-i}(\cO_{X}(*Y))\otimes\rd_{t}^{i},
\leqno(3.1.2)
$$
and similarly for
$ F_{p}\cB_{h} $.
(Here the shift of the index by
$ 1 $ associated to the direct image by a closed embedding of
codimension
$ 1 $ is omitted to simplify the notation.)

Let
$ V $ be the filtration of Kashiwara [22] and Malgrange [27] on
$ \cB_{h} $,
$ \cB_{h}[t^{-1}] $ along
$ X\times\{0\} $ indexed by
$ \bQ $.
This is an exhaustive decreasing filtration of coherent
$ \cD_{X} $-submodules, and satisfies the following conditions:

\medskip\noindent
(i)
$ t(V^{\alpha}\cB_{h}[t^{-1}]) \subset V^{\alpha +1}\cB_{h}[t^{-1}] $,
$ \rd_{t}(V^{\alpha}\cB_{h}[t^{-1}]) \subset V^{\alpha -1}
\cB_{h}[t^{-1}] $ for
$ \alpha \in \bQ $,

\medskip\noindent
(ii)
$ t(V^{\alpha}\cB_{h}[t^{-1}]) = V^{\alpha +1}\cB_{h}[t^{-1}] $ for
$ \alpha > 0 $,

\medskip\noindent
(iii)
$ \rd_{t}t - \alpha $ is nilpotent on
$ \Gr_{V}^{\alpha}\cB_{h}[t^{-1}] $ for
$ \alpha \in \bQ $,

\medskip\noindent
and similarly for
$ \cB_{h} $.
Here
$ \Gr_{V}^{\alpha} = V^{\alpha}/V^{>\alpha} $ with
$ V^{>\alpha} = \bigcup_{\beta>\alpha}V^{\beta} $, and we assume
$ V $ is indexed discretely and satisfies
$ V^{\alpha} = V^{\alpha-\varepsilon} $ for
$ \varepsilon > 0 $ sufficiently small.
Note that the filtration
$ V $ on
$ \cB_{h} $ is induced by that of
$ \cB_{h}[t^{-1}] $, and
$$
(V^{\alpha}\cB_{h},F) = (V^{\alpha}\cB_{h}[t^{-1}],F)\quad
\hbox{for}\,\,\,\alpha > 0,
\leqno(3.1.3)
$$
because
$ \cB_{h}[t^{-1}]/\cB_{h} $ is supported on
$ X\times\{0\} $ so that
$ \Gr_{V}^{\alpha}(\cB_{h}[t^{-1}]/\cB_{h}) = 0 $
unless
$ -\alpha \in \bN $.
We will denote also by
$ V $ the filtration on
$ \cO_{X}(*Y) $ induced by that on
$ \cB_{h}[t^{-1}] $.

\medskip\noindent
{\bf 3.2.~Proposition.} {\it
Let
$ j : X\times (\bA^{1}\setminus\{0\}) \to X\times\bA^{1} $
denote the inclusion.
Then
}
$$
F_{0}(\cB_{h}[t^{-1}]) =
V^{0}(\cB_{h}[t^{-1}])\cap j_{*}j^{*}F_{0}\cB_{h},
\leqno(3.2.1)
$$
$$
F_{0}(\cO_{X}(*Y)) = V^{0}(\cO_{X}(*Y)) = V^{0}(\cO_{X}(Y)).
\leqno(3.2.2)
$$

\medskip\noindent
{\it Proof.}
The isomorphism (3.2.1) follows from the property of the Hodge
filtration of a mixed Hodge module on which the action of
$ t $ is bijective (see [30], 4.2) because
$ \min\{p\in\bZ\,|\,F_{p}(\cB_{h}[t^{-1}])\ne 0\}=0 $.
Then (3.2.1) implies the first isomorphism of (3.2.2), because
$ F_{0}(\cB_{h}[t^{-1}])|_{U'} = V^{0}(\cO_{X}(*Y))|_{U'} =
\cO_{X}(Y)|_{U'} $ where
$ U' = X\setminus\Sing\,Y $.
The last isomorphism of (3.2.2) is equivalent to
$ V^{k}\cO_{X}\subset h^{k-1}\cO_{X} $ for any positive integer
$ k $ because
$ gh^{-k}\in V^{0}(\cO_{X}(*Y)) $ with
$ g\in \cO_{X} $ is equivalent to
$ g\in V^{k}\cO_{X} $.
So it is proved by restricting to the smooth points of
$ Y_{\rm red} $ because
$ \cO_{X,x} $ is a unique factorization domain.

\medskip\noindent
{\bf 3.3.~Proof of Theorem~2.}
We have
$ \tV^{>0}(\omega_{X}(Y)) \supset \omega_{X} $ because
$ \tV^{>1}\cO_{X} \supset V^{>1}\cO_{X} = \cO_{X}(-Y) $.
So the first isomorphisms in the first and the second rows
of Theorem~2 follow from the vanishing of
$ H^{k}(X,\omega_{X}) $ for
$ k = 0, 1 $.
Then the second isomorphism of the first row follows from
the last assertion of Proposition~3.2.
We now reduce the remaining assertions to
$$
\tV^{>0}\omega_{Y} = \tom_{Y} := \rho_{*}\omega_{\tY},\quad
\tV^{>1}\omega_{X} = \rho_{*}\omega_{\tX}(-E),
\leqno(3.3.1)
$$
where
$ \rho : \tX \to X $ is an embedded resolution of singularities of
$ Y $ such that
$ \rho^{-1}(Y) $ and the exceptional divisor
$ E $ are divisors with normal crossings and the proper transform
$ \tY $ of
$ Y $ is smooth.
Here the multiplicities of
$ E $ is defined by
$ \rho^{*}(Y) = E + \tY $, and the adjoint ideal is defined by
$ \rho_{*}\omega_{\tX/X}(-E) $.

Assuming (3.3.1) and using, for example, a cubic resolution,
the last isomorphism in the second row is more or less well-known.
The relation with the intersection cohomology is reduced to the
assertion that
$ \tom_{Y} $ gives the first nontrivial piece of the Hodge
filtration of the mixed Hodge module corresponding to the
intersection complex.
Here the last assertion follows easily from the decomposition
theorem in the category of mixed Hodge modules because
$ \tom_{Y} $ is torsion-free.

We now show (3.3.1).
By [17] we have the exactness of the short exact sequence
$$
0 \to \rho_{*}\omega_{\tX} \to \rho_{*}\omega_{\tX}(\tY) \to
\rho_{*}\omega_{\tY} \to 0.
$$
So the two assertions in (3.3.1) are equivalent to each other,
because the second assertion is equivalent to
$ \tV^{>0}(\omega_{X}(Y)) = \rho_{*}\omega_{\tX}(\tY) $, see also
[5], [14], [38].
We will show the first assertion of (3.3.1).

Since the assertion is local and we have a local equation
$ h $ of
$ Y $, we may consider
$ \tV^{>1}\omega_{X} $ instead of
$ \tV^{>0}(\omega_{X}(Y)) $ by trivializing
$ \cO_{X}(Y) $.
By definition
$ \tV^{>1}\omega_{X} $ contains
$ V^{>1}\omega_{X} = \omega_{X}(-Y) $ and
$$
\tV^{>1}\omega_{X}/V^{>1}\omega_{X} \subset
\Gr_{V}^{1}\omega_{X}
$$
is identified with
$$
\Ker(N:\Gr_{V}^{1}\cB_{h}\to\Gr_{V}^{1}\cB_{h})\otimes
\omega_{X}\cap\Gr_{V}^{1}\omega_{X},
$$
because
$ t : \Gr_{V}^{0}\cB_{h}\to\Gr_{V}^{1}\cB_{h} $ is injective.
(Here
$ N = t\rd_{t} $ as usual, and the tensor product with
$ \omega_{X} $ may be viewed as the transformation between
left and right
$ \cD $-modules.)
Then it is further identified with
$ \tom_{Y} \subset \Gr_{V}^{1}\omega_{X} $,
see the proof of Th. 0.6 in [30].
So the assertion follows.

\medskip\noindent
{\bf 3.4.~Proof of Theorem~3.}
The first isomorphism is shown in [32], 0.9, where the acyclicity
of the Koszul complex associated to the partial derivatives
$ \partial f/\partial x_{i} $ is used in an essential way.
The second isomorphism then follows from Theorem~1 by considering
the image of
$ \Gamma(X,F_{q}(\omega_{X}(*Y)))\subset
\Gamma(X,P_{q}(\omega_{X}(*Y))) $ by the last morphism of (0.1) in
the introduction.

\medskip\noindent
{\bf 3.5.~Remarks.} (i)
If
$ h $ is semi-weighted-homogeneous with weight
$ (w_{1},\dots,w_{n}) $ as in Introduction,
it is well known that
$ \alpha_{Y,y} = \sum_{i}w_{i} $.
If
$ h $ is weighted homogeneous, this is due to Kashiwara
(unpublished), and follows also from [26], [34] together with
[33], [39] (or we can use a calculation of the Gauss-Manin
connection by Brieskorn together with [26]).
If
$ h $ is semi-weighted-homogeneous, the assertion is then
reduced to the fact that
$ \alpha_{Y,y} $ does not change by a
$ \mu $-constant deformation in the isolated singularity
case.
The last assertion follows, for example, from the fact that
$ \alpha_{Y,y} $ coincides with the minimal
spectral number [35] by [26], [33], [39],
because the spectral numbers are constant
under a
$ \mu $-constant deformation [40].

\medskip
(ii) Let
$ \cG_{f}^{j} $ be the algebraic Gauss-Manin systems associated
to a homogeneous polynomial
$ f : \bA^{n+1}\to\bA^{1} $ as in the introduction.
These are the direct image sheaves of
$ \cB_{f} $ by the projection
$ \bA^{n+1}\times\bA^{1}\to\bA^{1} $, and are defined by using
the relative de Rham complex associated to the projection.
We identify these with the corresponding modules over the
Weyl algebra by taking the global section functor.
They have the filtration
$ V $, which satisfy the conditions similar to those in
(3.1), and are induced by the filtration
$ V $ on
$ \cB_{f} $ by using the relative de Rham complex.
Let
$ \delta(t-f) $ denote the canonical generator
$ 1\otimes 1 $ of
$ \cB_{f} $.
Let
$ b_{f}(s) $ be the Bernstein-Sato polynomial of
$ f $.
Let
$ \alpha_{f} $ and
$ \beta_{f} $ be respectively the minimal and the maximal root of
$ b_{f}(-s)/(-s+1) $.
Then
$ \alpha_{f} > 0 $ by [21],
$ \beta_{f} \le n + 1 - \alpha_{f} $ by [31], and
$$
V^{>\beta_{f}-1}\cG_{f}^{0}\subset\ocH_{f}\subset
V^{\alpha_{f}}\cG_{f}^{0}.
\leqno(3.5.1)
$$
This holds also for the analytic Brieskorn modules and the analytic
Gauss-Manin systems associated to a germ of a holomorphic
function
$ f $.
In the isolated singularity case, this is well-known, and
follows from [26] together with [33], [35], [39]
(see also [34] for the weighted homogeneous case).
In general, the first inclusion of (3.5.1) is reduced to
$$
\cD_{X}[s]\delta(t-f) \supset V^{>\beta_{f}-1}\cB_{f},
$$
and the latter follows from the fact that
$ b_{f}(s) $ is the minimal polynomial for the action of
$ s = -\rd_{t}t $ on
$$
\cD_{X}[s]\delta(t-f)/t\cD_{X}[s]\delta(t-f).
$$
The second inclusion of (3.5.1) follows from the fact that
the action of
$ \rd_{t} $ is bijective on
$ \cG_{f}^{0} $ so that the filtration
$ V $ can be replaced with the microlocal filtration in [31].

Note that the first inclusion of (3.5.1) is equivalent to
$$
\Gr_{V}^{\alpha}\ocH_{f} = \Gr_{V}^{\alpha}\cG_{f}^{0}
\quad\hbox{for any}\,\,\, \alpha > \beta_{f} - 1,
\leqno(3.5.2)
$$
and also to the isomorphisms
$$
t : \Gr_{V}^{\alpha}\ocH_{f} \simto \Gr_{V}^{\alpha+1}\ocH_{f}
\quad\hbox{for any}\,\,\, \alpha > \beta_{f} - 1.
\leqno(3.5.3)
$$
Note also that (1.6.2) is equivalent to
$$
\Gr_{V}^{\alpha}\ocH_{f} = \ocH_{f,k} \quad\hbox{for}\quad
\alpha=k/d.
\leqno(3.5.4)
$$

\bigskip\bigskip
\centerline{\bf 4. Kodaira-Spencer map}

\bigskip\noindent
{\bf 4.1.~Deformation of smooth open varieties.}
Let
$ \pi : \tX \to S $ be a proper smooth morphism of complex
manifolds, and
$ \tY $ be a divisor with normal crossings which is flat over
$ S $ and such that
$ \tY_{s} := \tY \cap \pi^{-1}(s) $ is a divisor with normal
crossings on
$ \tX_{s} := \pi^{-1}(s) $ for any
$ s \in S $.
Put
$ U = \tX \setminus \tY $, and
$ U_{s} = \tX_{s} \setminus \tY_{s} $.

Fix
$ s \in S $, and let
$ \Theta_{\tX_{s}}(\log \tY_{s}) $ be the sheaf of logarithmic
vector fields on
$ \tX_{s} $ along
$ \tY_{s} $ which is the dual of
$ \Omega_{\tX_{s}}^{1}(\log \tY_{s}) $ (i.e.
$ \xi \in \Theta_{\tX_{s}} $ belongs to
$ \Theta_{\tX_{s}}(\log \tY_{s}) $ if and only if
$ \xi\cI_{\tY_{s}} \subset \cI_{\tY_{s}} $ where
$ \cI_{\tY_{s}} $ is the (reduced) ideal sheaf of
$ \tY_{s} $).
There is the Kodaira-Spencer map
$$
T_{S,s} \to H^{1}(\tX_{s},\Theta_{\tX_{s}}(\log \tY_{s})),
\leqno(4.1.1)
$$
as in the classical case (where
$ \tY = \emptyset $), see e.g. [23], [37].
Furthermore, it is known that the
$ \cO_{S} $-linear part of the
Gauss-Manin connection
$$
\Gr_{F}\nabla_{\xi} :
\Gr_{F}^{p}H^{j}(U_{s},\bC) \to
\Gr_{F}^{p-1}H^{j}(U_{s},\bC)
\leqno(4.1.2)
$$
coincides up to a sign with the action of the image of
$ \xi $ by the Kodaira-Spencer map, see loc. cit.
In the classical case, this is due to Griffiths, see [18].

\medskip\noindent
{\bf 4.2.~Case of the complement of hypersurfaces.}
Let
$ X_{s} = \bP^{n} $, and let
$ \{Y_{s}\}_{s\in S} $ be an equisingular family of divisors
defined by polynomials
$ f_{s} $ which depend algebraically on
$ s \in S $, where
$ S $ is assumed to be a smooth affine variety.
Here equisingular means that
$ (X_{s},Y_{s})_{s\in S} $ admits a simultaneous embedded resolution
$ (\tX_{s},\tY_{s}) \to (X_{s},Y_{s}) \, (s\in S) $ which is
induced by an embedded resolution
$ (\tX,\tY) \to (X,Y) $ by restricting to the fibers at each
$ s \in S $, where
$ X = \bP^{n}\times S $ and
$ Y = \bigcup_{s\in S} Y_{s}\times\{s\} $.
This assumption implies that
$ \{U_{s}\}_{s\in S} $ is topologically locally trivial, and
$ \{H^{j}(U_{s},\bC)\}_{s\in S} $ (or, more precisely,
$ R^{j}\pi_{*}\bC_{U} $) is a local system on
$ S $ which underlies naturally a variation
of mixed Hodge structure, see [29].
In particular, the dimension of the Hodge filtration
$ F^{k}H^{n}(U_{s},\bC) $ for
$ s \in S $ is constant .

\medskip\noindent
{\bf 4.3.~Remarks.} (i)
It is not clear whether the dimension of the pole order filtration
$ P^{n-q} $ on
$ H^{n}(U_{s},\bC) $ is constant for
$ s \in S $.
By Theorem~1 this is equivalent to that
$ \dim\ocH_{f_{s},qd} $ is constant.
Note that
$ \dim (\rd f_{s})_{qd-n-1}\, (s\in S) $ is not necessarily
constant, where
$ (\rd f_{s})_{j} $ is the degree
$ j $ part
of the Jacobian ideal
$ (\rd f_{s}) $.
For example, consider a family of plane curves defined by
$ f_{s} = x^{4}z + y^{5} + sx^{2}y^{3} = 0 $ for
$ s \in \bC $.
Here
$ \Sing\,Y_{s} $ is one point, and the family is equisingular
in the sense of (4.2) (indeed, the resolution of singularities
in the case of irreducible plane curves depends only on the
Puiseux pairs).
Then the local Tjurina number jumps at
$ s = 0 $, and
$ \dim (\rd f_{s})_{k}\, (s\in S) $ is not constant at
$ s = 0 $ for
$ k \gg 0 $, see [7].

\medskip
(ii) There is a theory of versal family for deformations of varieties
with normal crossing divisors as in (4.1), see [23].
If
$ Y_{s} $ has only isolated singularities (where
$ s $ is a base point of
$ S $), then one can apply the above theory to a deformation
of a resolution of singularities
$ (\tY_{s},E_{s}) \to (Y_{s},\Sigma_{s}) $ (where
$ E_{s} $ is the exceptional divisor and
$ \Sigma_{s} = \Sing\,Y_{s} $)
instead of applying it to the embedded resolution
$ (\tX_{s},\tY_{s})\to (X_{s},Y_{s}) $ as in (4.2).
In certain cases, it is possible to blow down a deformation
$ (\tY_{s'},E_{s'}) $ of
$ (\tY_{s},E_{s}) $ by [16], e.g. if
$ \dim Y_{s} = 2 $ and the singularities of
$ Y_{s} $ are rational double points so that the exceptional
divisor
$ E_{s} $ is a disjoint union of copies of
$ \bP^{1} $.
This has an advantage that the dimension of
$ \tY_{s} $ is smaller than that of
$ \tX_{s} $.
However, it it not clear whether the blow-down of an arbitrary
deformation
$ (\tY_{s'},E_{s'}) $ of
$ (\tY_{s},E_{s}) $ is still a hypersurfaces of
$ \bP^{n} $, and we have to determine the subset of
the base space of the versal deformation of
$ \tY_{s} $ consisting of the points corresponding to
hypersurfaces of
$ \bP^{n} $.

There is also a problem about the difference between the
versal deformation of
$ (\tY_{s},E_{s}) $ and that of
$ \tY_{s} $.
Assume, for simplicity,
$ E_{s} $ is smooth, and let
$ N_{E_{s}} $ denote the normal bundle of
$ E_{s} $ in
$ \tY_{s} $.
Then there is a short exact sequence
$$
0\to \Theta_{\tY_{s}}(\log E_{s})\to \Theta_{\tY_{s}}\to
N_{E_{s}}\to 0,
\leqno(4.3.1)
$$
inducing a long exact sequence
$$
\eqalign{
H^{0}(E_{s},N_{E_{s}})
&\to H^{1}(\tY_{s},\Theta_{\tY_{s}}(\log E_{s}))
\to H^{1}(\tY_{s},\Theta_{\tY_{s}})\cr
\buildrel{\gamma}\over\to H^{1}(E_{s},N_{E_{s}})
&\to H^{2}(\tY_{s},\Theta_{\tY_{s}}(\log E_{s}))
\to H^{2}(\tY_{s},\Theta_{\tY_{s}}).\cr
}
\leqno(4.3.2)
$$
This may be used to study the difference between the versal
deformation of
$ \tY_{s} $ and that
$ (\tY_{s},E_{s}) $ (the latter may be viewed as
the versal equisingular deformation of
$ (Y_{s},\Sigma_{s}) $ in the sense of (4.2) if we can blow
down as above).
We have usually
$ H^{0}(E_{s},N_{E_{s}}) = 0 $ (e.g. if
$ N_{E_{s}} $ is negative).
We have sometimes
$ H^{1}(E_{s},N_{E_{s}}) = 0 $ for example if
$ E_{s} = \bP^{m-1} $ with
$ m := \dim \tY_{s} > 2 $ by the Bott vanishing theorem.
However, the morphism
$ \gamma $ in (4.3.2) does not vanish in general.
For example, if
$ \omega_{\tY_{s}} $ is trivial, then
$ \gamma $ is identified with the restriction morphism
$$
H^{1}(\tY_{s},\Omega^{m-1}_{\tY_{s}})\to
H^{1}(E_{s},\Omega^{m-1}_{\tY_{s}}),
$$
and this morphism is surjective in the case
$ m = 2 $.
This may be related to [6].

\medskip\noindent
{\bf 4.4.~Gauss-Manin systems.}
Let
$ Z = \bA^{n+1} $, and
$ f $ be a function on
$ Z\times S $ whose restriction to
$ Z\times\{s\} $ is
$ f_{s} $.
Let
$ \cG_{(f,pr)} $ be the direct image of
$ \cO_{Z\times S} $ by
$ (f,pr) : Z\times S \to \bA^{1}\times S $ as an algebraic
$ \cD $-module.
This is obtained by taking the global section functor of the
relative de Rham complex
$ \DR_{Z\times S/S}(\cB_{f}) $.
Here we may take the global section functor because
$ \bA^{1}\times S $ is affine.
Recall that
$$
\cB_{f} = \cO_{Z\times S}\otimes_{\bC}\bC[\rd_{t}],
$$
and the actions of
$ \xi \in \Theta_{Z\times S} $ and
$ t $ are given by
$$
\eqalign{
\xi(g\otimes\rd^{i})
&= (\xi g)\otimes\rd^{i} - (\xi f)g\otimes\rd^{i+1},\cr
t(g\otimes\rd^{i})
&= fg\otimes\rd^{i} - ig\otimes\rd^{i-1},\cr
}
\leqno(4.4.1)
$$
where the actions of
$ \cO_{Z\times S} $ and
$ \rd_{t} $ are natural ones.

Let
$ \cG_{f_{s}} $ be the restriction of
$ \cG_{(f,pr)} $ to
$ \bA^{1}\times \{s\} $ which is by definition the tensor product of
$ \cG_{(f,pr)} $ with
$ \Gamma(\bA^{1}\times \{s\},\cO) $ over
$ \Gamma(\bA^{1}\times S,\cO) $.
By the assumption on the equisingularity,
$ \cG_{f_{s}} $ is the direct image of
$ \cB_{f_{s}} $ which is defined by replacing
$ f $ with
$ f_{s} $ and
$ S $ with
$ \{s\} $.
There are decompositions
$$
\cG_{(f,pr)} = \mopls_{k\in\bZ}\cG_{(f,pr),k}, \quad
\cG_{f_{s}} = \mopls_{k\in\bZ}\cG_{f_{s},k},
$$
such that the degrees of
$ x_{i} $ and
$ dx_{i} $ are
$ 1 $, where the
$ x_{i} $ are the coordinates of
$ Z := \bA^{n+1} $.

We see that
$ \cG_{f_{s},k} $ is the tensor product of
$ \cG_{(f,pr),k} $ and
$ \Gamma(\bA^{1}\times \{s\},\cO) $ over
$ \Gamma(\bA^{1}\times S,\cO) $, and the canonical morphism
$$
\cH_{f_{s}}\to\cG_{f_{s}}
\leqno(4.4.2)
$$
is compatible with the decomposition.
Note that the action of
$ \nabla_{t\rd/\rd t} $ on
$ \cG_{f_{s},k} $ is the multiplication by
$ k/d - 1 $, and the image of (4.4.2) is
$ \ocH_{f_{s}} $, see [1].
By the assumption on the equisingularity,
$ \dim \cG_{f_{s},k} $ is constant, and hence
$ \cG_{(f,pr),k} $ is a projective
$ R $-module, where
$ R = \Gamma(S,\cO_{S}) $.
However, this is not clear for
$ \ocH_{f_{s},k} $, see Remark (4.3)(i).

\medskip\noindent
{\bf 4.5.~Theorem.} {\it Assume
$ \dim P^{j}H^{n}(U_{s},\bC) \,(s\in S) $ is constant for any
$ j \in \bZ $.
Then for
$ \xi\in\Theta_{S} $, the action of the Gauss-Manin
connection
$ \nabla_{\xi} $ induces a well-defined morphism
$$
\Gr_{P}\nabla_{\xi} : \Gr_{P}^{n-q+1}H^{n}(U_{s},\bC) \to
\Gr_{P}^{n-q}H^{n}(U_{s},\bC),
\leqno(4.5.1)
$$
which corresponds to the multiplication by
$ -q(\xi f)_{s} $ using the identification in Theorem~{\rm 1}.
Furthermore there is a commutative diagram
$$
\matrix{
\Gr_{F}^{n-q+1}H^{n}(U_{s},\bC) &
\buildrel{{\rm Gr}_{F}\nabla_{\xi}}\over\longrightarrow &
\Gr_{F}^{n-q}H^{n}(U_{s},\bC)\cr
\downarrow && \downarrow \cr
\Gr_{P}^{n-q+1}H^{n}(U_{s},\bC) &
\buildrel{{\rm Gr}_{P}\nabla_{\xi}}\over\longrightarrow &
\Gr_{P}^{n-q}H^{n}(U_{s},\bC)\cr
\Vert && \Vert \cr
\ocH_{f,qd}/f\ocH_{f,(q-1)d} &
\buildrel{-q(\xi f)_{s}}\over\longrightarrow &
\ocH_{f,(q+1)d}/f\ocH_{f,qd}\cr
}
$$
where the vertical morphisms are either the natural morphisms or
the isomorphisms given by Theorem~{\rm 1}.
}

\medskip\noindent
{\it Proof.}
Let
$ \ocH_{(f,pr),k} $ be the image of the canonical morphism
$$
(\bC[x]_{k-n-1}\otimes_{\bC} R\otimes 1)
dx_{0}\wdg\cdots\wdg dx_{n} \to \cG_{(f,pr),k},
$$
where
$ \bC[x]_{j} $ is the degree
$ j $-part of the polynomial ring
$ \bC[x] $.
Since the tensor product is right exact, the
$ \ocH_{f_{s},k} $ are given by the image of the tensor product of
the inclusion
$$
\ocH_{(f,pr),k} \to \cG_{(f,pr),k}
$$
with
$ \bC_{s} := R/m_{s} $ over
$ R $ where
$ m_{s} $ is the maximal ideal at
$ s $.
Since
$ \dim P^{j}H^{n}(U_{s},\bC) $ is constant, we see that
$ \ocH_{(f,pr),k} $ is a projective
$ R $-module.
Let
$$
\ocH_{(f,pr)} = \mopls_{k\in\bN} \ocH_{(f,pr),k}\subset
\cG_{(f,pr)}.
$$
Then
$ \ocH_{(f,pr)} $ is stable by the action of
$ t\xi $ for
$ \xi\in\Theta_{S} $.
Indeed, we have by (4.4.1)
$$
\eqalign{
t\xi(g\otimes 1)
&= t((\xi g)\otimes 1 - (\xi f)g\otimes \rd_{t})\cr
&= f(\xi g)\otimes 1 - t\rd_{t}((\xi f)g\otimes 1),\cr
}
$$
and the action of
$ t\rd_{t} $ on
$ ((\xi f)g\otimes 1)dx_{0}\wdg\cdots\wdg dx_{n} $ is the
multiplication by
$ k/d $ if
$ g $ has pure degree
$ k - n - 1 $.
Calculating
$ \xi(gf^{-q}dx_{0}\wdg\cdots\wdg dx_{n}) $ we see that
the action of
$ \xi $ is compatible with the isomorphism in Theorem~1.
Since the action of the Gauss-Manin connection is
compatible with the structure of
$ \cD $-modules, we get the well-definedness of
$ \Gr_{P}\nabla_{\xi} $ together with the above commutative
diagram.
This completes the proof of (4.5).

\medskip\noindent
{\bf 4.6.~Remark.}
In general it is quite difficult to get a good condition for
the constantness of
$ \dim P^{j}H^{n}(U_{s},\bC) \,(s\in S) $ for any
$ j \in \bZ $.
If
$ F^{j} = P^{j} $ on
$ H^{n}(U_{s},\bC) $ for
$ j = k $,
$ k+1 $ and
$ \dim P^{k-1}H^{n}(U_{s},\bC) $ is constant (e.g. if
$ \dim Y_{s} = 2 $,
$ d = 4 $,
$ k = 2 $ and
$ P^{1}H^{2}(U_{s},\bC) = H^{2}(U_{s},\bC) $), then we can
calculate by (4.5)
$$
\Gr_{F}\nabla_{\xi} : \Gr_{F}^{k}H^{n}(U_{s},\bC) \to
\Gr_{F}^{k-1}H^{n}(U_{s},\bC) \subset
\Gr_{P}^{k-1}H^{n}(U_{s},\bC).
$$

\vfill\eject
\centerline{{\bf References}}

\bigskip
{\mfont
\item{[1]}
D.~Barlet and M.~Saito, Brieskorn modules and Gauss-Manin systems
for non isolated hypersurface singularities, preprint
(math.CV/0411406).

\item{[2]}
A.~Beilinson, J.~Bernstein and P.~Deligne, Faisceaux Pervers,
Ast\'erisque, vol.~100, Soc. Math. France, Paris, 1982.

\item{[3]}
E.~Brieskorn, Die Aufl\"osung der rationalen Singularit\"aten
holomorpher Abbildungen, Math. Ann. 178 (1968), 255--270.
\item{[4]}
E.~Brieskorn, Die Monodromie der isolierten Singularit\"aten von
Hyperfl\"achen, Manuscripta Math., 2 (1970), 103--161.

\item{[5]}
N.~Budur and M.~Saito, Multiplier ideals,
$ V $-filtration, and spectrum,
J. Algebraic Geom. 14 (2005), 269--282.

\item{[6]}
D.M.~Burns and J.M.~Wahl,
Local contributions to global deformations of surfaces,
Invent. Math. 26 (1974), 67--88.

\item{[7]}
A.D.R.~Choudary and A.~Dimca,
Koszul Complexes and hypersurface singularities,
Proc. AMS 121 (1994), 1009-1016.

\item{[8]}
P.~Deligne,
Equations Diff\'erentielles \`a Points Singuliers R\'eguliers,
Lect. Notes in Math. vol.~163, Springer, Berlin, 1970.

\item{[9]}
P.~Deligne, Th\'eorie de Hodge I, Actes Congr\`es
Intern. Math., 1970, vol. 1, 425-430; II, Publ. Math. IHES,
40 (1971), 5--57; III, ibid., 44 (1974), 5--77.

\item{[10]}
P.~Deligne and A.~Dimca,
Filtrations de Hodge et par l'ordre du p\^ole pour les hypersurfaces
singuli\`eres, Ann. Sci. Ecole Norm. Sup. (4) 23 (1990), 645--656.

\item{[11]}
A.~Dimca,
On the Milnor fibrations of weighted homogeneous polynomials,
Compositio Math. 76 (1990),19--47.

\item{[12]}
A.~Dimca,
Singularities and Topology of Hypersurfaces, Springer 1992.

\item{[13]}
A.~Dimca, Hodge numbers of hypersurfaces,
Abh. Math. Sem. Univ. Hamburg 66 (1996), 377--386.

\item{[14]}
L.~Ein and R.~Lazarsfeld, Singularities of theta divisors
and the birational geometry of irregular varieties,
J. Amer. Math. Soc. 10 (1997), 243--258.

\item{[15]}
M.~Goresky and R. MacPherson, Intersection homology theory,
Topology 19 (1980), 135--162.

\item{[16]}
H.~Grauert, \"Uber Modifikationen und exzeptionelle analytische
Mengen, Math. Ann. 146 (1962), 331--368.

\item{[17]}
H.~Grauert and O.~Riemenschneider, Verschwindungss\"atze
f\"ur analytische Kohomologiegruppen auf Komplexen R\"aumen,
Inv. Math. 11 (1970), 263--292.

\item{[18]}
M.~Green, Infinitesimal method in Hodge theory, in Algebraic cycles
and Hodge theory, Lect. Notes in Math. vol. 1594, Springer,
Berlin, 1994, pp. 1--92.

\item{[19]}
P.~Griffiths, On the period of certain rational integrals I,
II, Ann. Math. 90 (1969), 460--541.

\item{[20]}
H.A.~Hamm, Z\"ur analytischen und algebraichen Beschreibung der
Picard-Lefchetz Monodromie, Habilitationsschrift, G\"ottingen, 1974.

\item{[21]}
M.~Kashiwara,
$ B $-functions and holonomic systems, Inv. Math. 38 (1976/77),
33--53.

\item{[22]}
M.~Kashiwara, Vanishing cycle sheaves and holonomic systems
of differential equations, Algebraic geometry (Tokyo/Kyoto,
1982), Lect. Notes in Math. 1016, Springer, Berlin,
1983, pp. 134--142.

\item{[23]}
Y.~Kawamata, On deformations of compactifiable complex manifolds,
Math. Ann. 235 (1978), 247--265.

\item{[24]}
R.~Lazarsfeld, Positivity in algebraic geometry II,
Springer, Berlin, 2004.

\item{[25]}
F.~Loeser, Quelques cons\'equences locales de la th\'eorie de
Hodge, Ann. Inst. Fourier 35 (1985), 75--92.

\item{[26]}
B.~Malgrange, Le polyn\^ome de Bernstein d'une
singularit\'e isol\'ee, in Lect. Notes in Math. 459, Springer,
Berlin, 1975, pp. 98--119.

\item{[27]}
B.~Malgrange, Polyn\^ome de Bernstein-Sato et cohomologie
\'evanescente, Analysis and topology on singular spaces, II,
III (Luminy, 1981), Ast\'erisque 101--102 (1983), 243--267.

\item{[28]}
M.~Merle and B.~Teissier,
Conditions d'adjonction d'apr\`es Du Val,
in: S\'eminaire sur les Singularit\'es des Surfaces, Springer
Lecture Notes in Math. 777 (1980), pp. 229--245.

\item{[29]}
M.~Saito, Mixed Hodge modules, Publ. RIMS, Kyoto Univ. 26
(1990), 221--333.

\item{[30]}
M.~Saito, On
$ b $-function, spectrum and rational singularity,
Math. Ann. 295 (1993), 51--74.

\item{[31]}
M.~Saito, On microlocal
$ b $-function. Bull, Soc. Math. France 122 (1994), 163--184.

\item{[32]}
M.~Saito, On the Hodge filtration of Hodge modules,
RIMS-preprint 1078, May 1996.

\item{[33]}
J.~Scherk and J.H.M.~Steenbrink,
On the mixed Hodge structure on the cohomology of the Milnor
fiber, Math. Ann. 271 (1985), 641--665.

\item{[34]}
J.H.M.~Steenbrink,
Intersection form for quasi-homogeneous singularities,
Compos. Math. 34 (1977), 211--223.

\item{[35]}
J.H.M.~Steenbrink, Mixed Hodge structure on the vanishing
cohomology, in Real and Complex Singularities (Proc. Nordic
Summer School, Oslo, 1976) Alphen a/d Rijn: Sijthoff-Noordhoff
1977, pp. 525--563.

\item{[36]}
J.H.M.~Steenbrink, Adjunction conditions for 1-forms on surfaces
in projective three-space, preprint (math.AG/0411405).

\item{[37]}
S.~Usui, Variation of mixed Hodge structures arising from family of
logarithmic deformations, Ann. Sci. Ecole Norm. Sup. (4) 16 (1983),
91--107.

\item{[38]}
M.~Vaqui\'e, Irr\'egularit\'e des rev\^etements cycliques,
in Singularities (Lille, 1991), London Math. Soc. Lecture
Note Ser., 201, Cambridge Univ. Press, Cambridge, 1994,
pp. 383--419.

\item{[39]}
A.~Varchenko, The asymptotics of holomorphic forms
determine a mixed Hodge structure, Soviet Math. Dokl. 22 (1980),
772--775.

\item{[40]}
A.~Varchenko,
The complex singularity index does not change along the stratum
$ \mu = $ const, Funk. Anal. 16 (1982), 1--12.

\bigskip
Alexandru Dimca

Laboratoire J.A. Dieudonn\'e, UMR du CNRS 6621, Math\'ematiques

Universit\'e de Nice-Sophia Antipolis

Parc Valrose, 06108 Nice Cedex 02, FRANCE

e-mail: dimca@math.unice.fr

\medskip
Morihiko Saito

RIMS Kyoto University, Kyoto 606--8502 JAPAN

e-mail: msaito@kurims.kyoto-u.ac.jp
}
\bye